\documentclass{elsart}
\usepackage{amssymb}
\usepackage{amsfonts}
\newcommand{\trace}{\mathop{\rm Tr}\nolimits}
\newcommand{\twomat}[4]{\left(\begin{array}{cc}{#1}&{#2}\\{#3}&{#4}\end{array}\right)}

\newcommand{\R}{{\mathbb{R}}}
\newcommand{\M}{{\mathbb{M}}}
\DeclareRobustCommand\openone{\leavevmode\hbox{\small1\normalsize\kern-.33em1}}

\newcommand{\be}{\begin{equation}}
\newcommand{\ee}{\end{equation}}
\newcommand{\bea}{\begin{eqnarray}}
\newcommand{\eea}{\end{eqnarray}}
\newcommand{\beas}{\begin{eqnarray*}}
\newcommand{\eeas}{\end{eqnarray*}}

\newtheorem{theorem}{Theorem}
\newtheorem{lemma}{Lemma}

\begin{document}
\begin{frontmatter}
\title{Spectral radius of Hadamard product versus conventional product for non-negative matrices}
\author{Koenraad M.R.\ Audenaert}
\address{
Dept.\ of Mathematics,\\
Royal Holloway, University of London,\\
Egham TW20 0EX, United Kingdom}
\ead{koenraad.audenaert@rhul.ac.uk}
\date{\today}
\begin{abstract}
We prove an inequality for the spectral radius of products of non-negative matrices conjectured by X.~Zhan. 
We show that for all $n\times n$ non-negative matrices $A$ and $B$, 
$\rho(A\circ B)\le\rho((A\circ A)(B\circ B))^{1/2}\le\rho(AB)$, where $\circ$ represents the Hadamard product.
\end{abstract}

\end{frontmatter}
We denote by $A_{ij}$ the entry of a matrix $A$ in position $(i,j)$, and
by $\circ$ the Hadamard product, that is, the entrywise product $(A\circ B)_{ij} = A_{ij} B_{ij}$.
We shall be considered with entrywise positive and entrywise non-negative matrices, which we'll simply call
positive and non-negative matrices, and denote by $A>0$ and $A\ge0$, respectively.

The spectral radius $\rho(A)$ of a matrix $A$ is the largest modulus of the eigenvalues of $A$.
For non-normal matrices, the spectrum does not build norms, for the simple reason that a matrix can be unbounded
even when all its eigenvalues are zero. The Jordan block $\twomat{0}{1}{0}{0}$ serves as the prime example here.
In particular, the spectral radius is not a matrix norm. This shows up for example in the fact that the
spectral radius is not submultiplicative:
$\rho(AB)\le \rho(A)\rho(B)$ does not hold in general, not even when restricting to non-negative matrices.
A counterexample is given by the pair \cite{HJII}
$$
A=\twomat{0}{0}{1}{0}, B=\twomat{1}{1}{0}{1}.
$$
On the other hand,
for non-negative $A$ and $B$, the spectral radius is submultiplicative w.r.t.\ the Hadamard product:
$\rho(A\circ B) \le \rho(A)\rho(B)$ (\cite{HJII}, observation 5.7.4).

In the light of the few relevant statements that can be made about the spectral radius, 
X.~Zhan made the remarkable discovery that for pairs of square non-negative matrices, the spectral radius
of the Hadamard product is always bounded above by the spectral radius of the (conventional) matrix product.
This was presented in the talk \cite{zhan} and posed as a conjecture.

The purpose of this note is to prove this conjecture. In fact, we will prove a little more:
\begin{theorem}
For $n\times n$ non-negative matrices $A$ and $B$, 
$$
\rho(A\circ B)\le\rho((A\circ A)(B\circ B))^{1/2}\le\rho(AB).
$$
\end{theorem}

Before going into the proof, we note that there is no reasonable lower bound on $\rho(A\circ B)$ in terms of
$\rho(AB)$.
The pair $A$ and $B$ given above also shows that $\rho(A\circ B)$ can be zero while $\rho(AB)$ is not.
In fact, this can happen even when $A$ and $B$ both have non-zero spectral radius.
Take 
$$
A = \twomat{1}{x}{0}{1}, B=\twomat{1}{0}{y}{1}.
$$
Then $A$, $B$ and $A\circ B$ all have spectral radius 1, while $AB=\twomat{1+xy}{x}{y}{1}$.
Since $\det(AB)=1$ and $\trace(AB)=2+xy$, the spectral radius of $AB$ is unbounded.

\bigskip

Let us now present the promised proof of Theorem 1.
The proof is based on the following representation 
of the spectral radius of positive matrices:
\begin{lemma}
For $A\in \M_n$ such that $A>0$, $\rho(A) = \lim_{m\to\infty} (\trace A^m)^{1/m}$.
\end{lemma}
\textit{Proof.}
By a basic theorem in Perron-Frobenius theory (see e.g.\ Theorem 8.2.8 in \cite{HJI}), for any $A>0$ we have
$$
\lim_{m\to\infty} \left(\frac{A}{\rho(A)}\right)^m = xy^T,
$$
where $x$ and $y$ are certain positive vectors for which $x^T y=1$.
Taking the trace yields
$$
\lim_{m\to\infty} \trace\left(\frac{A}{\rho(A)}\right)^m = 1,
$$
and the lemma follows immediately.
\qed

The lemma is not generally true for non-negative matrices, because the limit need not even exist.
The simplest example is the matrix 
$$A=\twomat{0}{1}{1}{0},$$ 
whose spectral radius is 1.
While its even powers are equal to the identity matrix
and have trace equal to 2, its odd powers are equal to $A$ itself and have trace equal to 0.
More generally, every permutation matrix shows this kind of behaviour, with a periodicity given by the
least common multiple of all the cycle lengths of the corresponding permutation.
In these cases, replacing the limit by the limit superior would salvage the lemma, but we do not
know whether this will be true in general.
On the other hand, the statement of the lemma does still hold unmodified for primitive non-negative matrices.

However, for our purposes we do not need all these generalisations.
In fact, we will only explicitly prove $\rho(A\circ B)\le\rho(AB)$ for positive $A$ and $B$.
Because of the continuity of the spectral radius, the Theorem then 
follows for non-negative $A$ and $B$ as well. 

An obvious consequence of the lemma is that for $A>0$, and any positive integer $k$,
$\rho(A^k) = \rho(A)^k$. By continuity of the spectral radius, this is also true for $A\ge0$.

\textit{Proof of Theorem 1 for positive $A$ and $B$.}
We will prove that the following inequalities hold for any positive integer $k$:
$$
\trace((A\circ B)^{2k}) \le \trace((A\circ A)(B\circ B))^k \le \trace((AB)^{2k}).
$$
By taking the $(2k)$th root and taking the limit $k\to\infty$, 
this implies the statement of the theorem, by the lemma.

The left-hand side $\trace((A\circ B)^{2k})$ can be written as a $2k$-fold sum:
\beas
\trace((A\circ B)^{2k}) &=
\sum_{i_1,i_2,\ldots,i_{2k}}&
(A_{i_1i_2}B_{i_1i_2})\,\,
(A_{i_2i_3}B_{i_2i_3})
\ldots
(A_{i_{2k}i_1}B_{i_{2k}i_1}) \\
&=
\sum_{i_1,i_2,\ldots,i_{2k}} &
\left(A_{i_1i_2}B_{i_2i_3}
\ldots
A_{i_{2k-1}i_{2k}}B_{i_{2k}i_1}\right) \\
&\quad\quad\quad\times&
\left(B_{i_1i_2}A_{i_2i_3}
\ldots
B_{i_{2k-1}i_{2k}}A_{i_{2k}i_1}\right).
\eeas
Note that the alternation of the $A$ and $B$ factors in the last line is intentional.

The last expression can be seen as an inner product between two vectors in $\R_+^{n^{2k}}$, one with entries
$A_{i_1i_2}B_{i_2i_3}
\ldots
A_{i_{2k-1}i_{2k}}B_{i_{2k}i_1}$ and the other with entries $B_{i_1i_2}A_{i_2i_3}
\ldots
B_{i_{2k-1}i_{2k}}A_{i_{2k}i_1}$.
One sees that these two vectors are the same up to a permutation of the entries (as can be seen by
performing a cyclic permutation on the indices $i_1,i_2,\ldots,i_{2k}$). Thus, in particular, both vectors have
the same Euclidean norm.
Applying the Cauchy-Schwarz inequality, $\langle x,y\rangle \le \sqrt{\langle x,x\rangle\,\langle y,y\rangle}$,
then gives
\bea
\trace((A\circ B)^{2k}) &\le
\sum_{i_1,i_2,\ldots,i_{2k}}&
\left(A_{i_1i_2}B_{i_2i_3}
\ldots
A_{i_{2k-1}i_{2k}}B_{i_{2k}i_1}\right) \nonumber\\
&\quad\quad\quad\times&
\left(A_{i_1i_2}B_{i_2i_3}
\ldots
A_{i_{2k-1}i_{2k}}B_{i_{2k}i_1}\right).\label{eq:star}
\eea
Noting that the last line can be succinctly written as $\trace((A\circ A)(B\circ B))^k$ proves the first
inequality of the theorem.

Now consider $\trace((AB)^{2k})$. This can be written as a $4k$-fold summation:
\beas
\trace((AB)^{2k}) &=
\sum_{i_1,i_2,\ldots,i_{2k} \atop j_1,j_2,\ldots,j_{2k}}&
(A_{i_1i_2}B_{i_2i_3}
\ldots
A_{i_{2k-1}i_{2k}}B_{i_{2k}j_1}) \\
&\quad\quad\quad\times&(A_{j_1j_2}B_{j_2j_3}
\ldots
A_{j_{2k-1}j_{2k}}B_{j_{2k}i_1}).
\eeas
The crucial observation is that if we take all terms of this summation for which $i_1=j_1$, $i_2=j_2$,\ldots
$i_{2k}=j_{2k}$, then we obtain the right-hand side of (\ref{eq:star}).
The terms we left out are of course all positive, because all matrix elements of $A$ and $B$ are positive.
Therefore, we find that $\trace((AB)^{2k})$ is an upper bound on (\ref{eq:star}), and thus on
$\trace((A\circ B)^{2k})$. This proves the second inequality.
\qed

\end{document}